# REJOINDER: CONDITIONAL GROWTH CHARTS

By Ying Wei and Xuming He

*Columbia University and University of Illinois*

First of all, we would like to thank all the discussants for their encouragement of and insightful comments on our work. We are especially appreciative of their unique perspectives brought to the topic of conditional growth charts. The discussants raised a number of interesting questions from statistical as well as clinical points of view, some of which are broader and deeper than what we will be able to address fully in this rejoinder. Our rejoinder will focus on the issues most directly related to the statistical model, called the global model, used in our paper.

**1. Conditional, marginal or joint models.** Carroll and Ruppert correctly pointed out that semiparametric efficient estimation of marginal models in the context of quantile regression calls for further research. The current literature on semiparametric efficient estimation often relies on the Gaussian likelihood, so there is no direct analogue in the quantile model without a parametric likelihood. Consider estimating the median of a univariate distribution from a correlated sample $Y_1, \ldots, Y_n$. Even in this much simpler setting, it is unclear if we can find a uniformly more efficient estimator than the usual sample median. Efficiency bounds similar to those of Newey and Powell [3] are yet to be developed for longitudinal models.

We chose to use the conditional model mainly driven by the desire to take into account the subject's prior growth path. We can integrate out the prior growth path to revert to a marginal model, and if the within-subject correlation in the marginal model is indeed accounted for through the dependence on prior growth path, this might lead to efficient estimation of the marginal model, but we have not explored it in detail.

Thompson used simulation results to show that if a joint Gaussian model holds (possibly after transformations), then distribution-based estimates of quantiles are less variable than the quantile regression estimates, especially for $\tau$ near 0 or 1. This observation was also made in [2] in a simpler setting.

Here we have the classical bias-variance trade-off, so our preference depends on the available sample size. Distributional assumptions might be necessary when there are no sufficient data, but our empirical work on growth data (including but not limited to the Finnish growth data presented in the paper) has shown that joint normality is often unrealistic, so not only should we worry about bias, the uncertainty estimates from such assumptions also cannot be trusted.

**2. Regularization with splines.** In our work, we used regression splines with preselected knots. Not surprisingly, Carroll and Ruppert preferred penalized splines backed up by the work of Ruppert [4]. We agree with the general conclusions there, but would add one item to complement their points—knots are more intuitive than a single smoothing parameter. As reinforced by Pere's discussion, there is often good understanding about the growth pattern of children (and adults), which can be used to guide us in knot selection. This was the approach we took. In [6] cubic B-splines were used for the growth charts of height for the wider age group of 0 to 20. There, 16 uneven knots were placed to reflect growth patterns at infancy and at puberty. It is harder to use the human growth information in choosing a smoothing parameter in a penalized approach. Obviously, some subjectivity is both a blessing and a curse in nonparametric smoothing, so we cannot make any further claims. We took comfort in the fact that the estimated quantile functions in our example were not very sensitive to the exact placement of the knots.

**3. Questions on our illustrative subject.** In the paper we chose a particular child, Subject 646 in the Finnish growth data, to illustrate how conditional growth charts can be used for weight screening. Observing that this boy was somewhat below the earlier growth path at 0.46 year, but moved above the path at 0.61 year, Pere suspected a measurement error at 0.46 year. Without a measurement error, Pere commented, the higher reading at 0.61 year might simply reflect natural variation in the weight growth of this child, not any indication of abnormal weight growth at that stage. These comments highlighted the difficulty in weight screening, because when a child starts to show some sign of weight increase beyond his earlier growth path, one may have multiple explanations for what might be happening.

In our example, the weight of the child was at the 43rd percentile on the marginal growth chart at 0.37 year, but dropped toward the 25th percentile at 0.46 year. If the weight measurement at 0.46 year were maintained at the 43rd percentile, the same global model would estimate the conditional weight (at 0.61 year) as 9.07 kg and 9.41 kg, respectively, for $\tau = 0.9$ and 0.97. Therefore, the observed weight at 0.61 year would still be above the 90th percentile on the conditional growth chart.

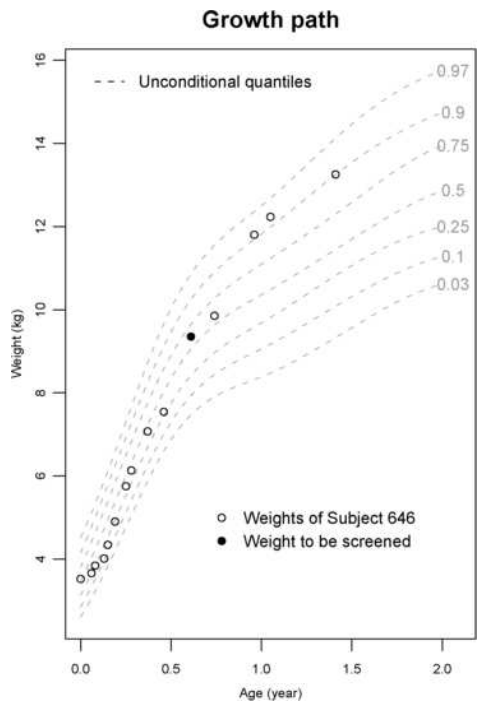

Fig. 1. *The full path of Subject* 646 *during infancy.*

To make an empirical assessment without the global model, we examined 96 subjects chosen as a reference group in the database. The subjects in this reference group had weight measurements between $0.37 - 0.04$ and $0.37 + 0.04$ year, and their weights were within 0.25 kg of Subject 646's weight at 0.37 year. At 0.61 year, the weight of Subject 646 was slightly above the 90th percentile of this group, and all the heavier boys were taller than Subject 646. One subject was taller by 1 cm, but all the others were taller by at least 2.5 cm. This analysis also confirms that by age 0.61 year the weight increase of Subject 646 should be taken as a flag even if the recorded weight at age 0.46 might have erred. In fact, we do have the actual weight measurements on the child at later ages. As shown in Figure 1, the child continued to move up in weight; after 1 year of age, his weight jumped above the 90th percentile curve in the marginal growth chart.

Several discussants questioned the use of height as a linear covariate in our conditional model. Clearly, we had similar concerns. The wide use of body mass index (BMI, defined as weight relative to height squared) suggested that height squared, $H^2$, might be more appropriate than height itself. Our empirical investigation suggested that the BMI is a useful index for adults, but not for young children. A similar statement was made in [1]. When we included both $H$ and $H^2$ in our global models for young children, the latter turned out to be statistically insignificant. Given the prior weight, it appeared that including the height linearly was quite appropriate. Our attempt to model the weight and height on the log scale did not lead to any better model.

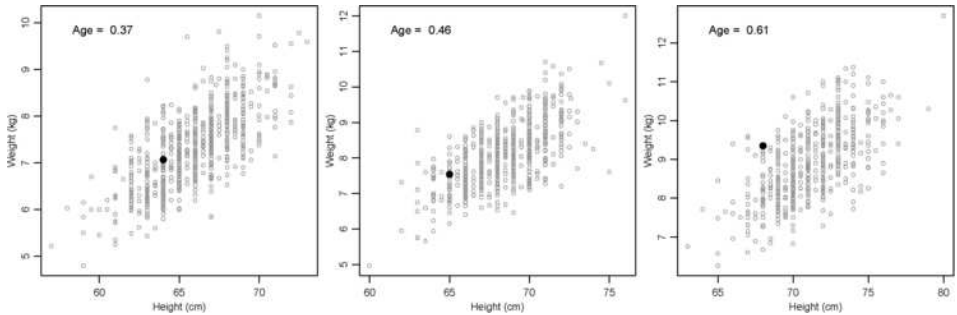

Fig. 2. *The solid dots are the actual height and weight measurements of Subject* 646 *at ages* 0.37, 0.46 *and* 0.61 *year.*

**4. Bivariate models.** Every discussant expressed preference to use height and weight jointly for screening, to which we fully concur. In fact, our conditional model was only a partial solution in this direction, and fully bivariate (or multivariate) solutions are clearly worth exploring. We appreciate Pere's description of [5], and thank Salibian-Barrera and Zamar for their extension to multivariate modeling. The connection with multivariate data depth is attractive, although it is unclear at this point which multivariate version of depth might be most appropriate. An even more interesting extension would be multivariate longitudinal assessment of growth, and the conditional models need to be further developed in this direction.

Figure 2 shows how Subject 646 compared to his peers in both weight and height at three time points. At each time point, boys with measurements within a window of 0.05 year were used as reference. It is clear from Figure 2 that this subject's height was consistently around the first quartile, but his weight started to be outlying by age 0.61. This is another indication that this subject was worth singling out for attention.

**5. Limitations.** Both Thompson and Pere raised questions about the diagnostic effectiveness of conditional growth charts. They have interpreted "screening" more clinically than we did.

Thompson correctly pointed out that even if a variable (e.g., weight) is strongly associated with an outcome (e.g., obesity), it might not be a useful screening tool. It is difficult to tie growth charts-based screening to a specific health problem. The growth charts, including the conditional ones, may be used as one of the screening tools for pediatricians, but we believe that any clinical diagnostic has to rely on additional information, and possibly medical tests.

We flagged out Subject 646 for his weight growth at 0.61 year in our example, and provided some further arguments for it in this rejoinder. However, this child, like all the others in the Finnish growth data, was considered to be normal for inclusion in the Finnish study. We can only observe that his weight growth started to fall out of the normal range, which would suggest to his parents and pediatrician to ask more questions. Possible changes in feeding or sleeping habit might be explored. Maybe the pediatrician could suggest closer follow-up of his weight to watch for possible signs of diabetes. Clearly, we do not have the expertise to make any specific clinical recommendation, and we do not expect growth charts alone to do the job. In this sense, diagnostic accuracy might have to be assessed in connection with additional tools available to the users of the growth charts.

Thompson's concern about once-off use of conditional charts was a real one, because the conditional model uses the subject's prior growth path as a starting point. In practice, growth charts are used "longitudinally." If the first sign of abnormal growth is observed at 12 months based on past measurements at 6 and 9 months, then the follow-up analysis at, say, 14 months, might use the measurements at 6 and 9 months as the covariates in the model. To this end, the flexibility of the global model to use a selected prior growth path is a plus.

Carroll and Ruppert suggested an alternative model to capture catch-up growth more specifically than our global model. Their model (3) is equivalent to modeling the rate of change,

$$(W_{ij} - W_{i,j-1})/D_{ij} = (g_\tau(t_{ij}) - g_\tau(t_{i,j-1}))/D_{ij} + b\{W_{i,j-1} - g_\tau(t_{i,j-1})\} + e_{i,j}.$$

The sign (and the size) of $b$ can be interpreted as catch-up or catch-down growth. This is such an intriguing idea that we explored it a little further. By taking $\tau = 0.5$ and $g_\tau$ to be the median in the marginal growth chart, we estimated $b(t)$ as shown in Figure 3. We see that the catch-up (corresponding to negative values of $b$) was mostly in the first three months. A similar model could be applied to $z$-scores of weight instead of weight itself, but the resulting estimate of $b(t)$ showed a similar pattern. On the other hand, we note that the catch-up growth may occur at different ages for different children, and the phenomenon is more evident for low-weight children during their infancy. In the model, $b(t)$ serves as a population summary. We also note that this model is less interpretable, and possibly inappropriate, at other $\tau$ values, and therefore it could not be used to estimate conditional quantiles for screening.

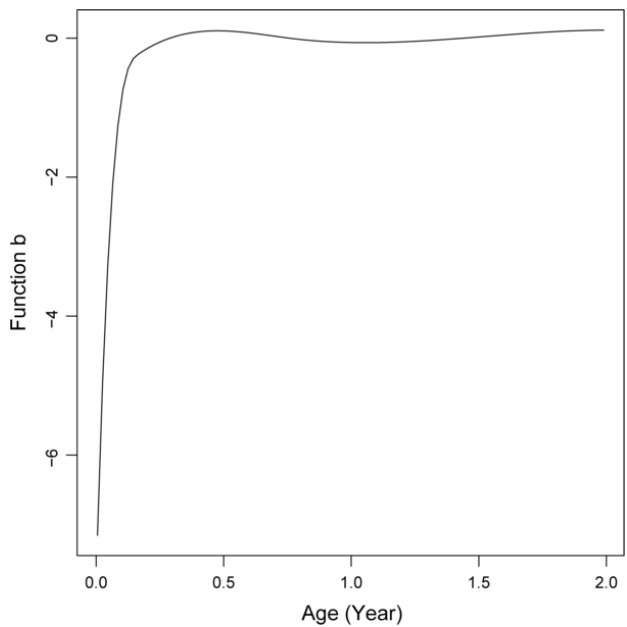

Fig. 3. *Estimated $b(t)$ as a cubic B-spline.*

## REFERENCES

[1] Franklin, M. F. (1999). Comparison of weight and height relations in boys from 4 countries. *American J. Clinical Nutrition* **70** 157–162.
[2] He, X. (1997). Quantile curves without crossing. *Amer. Statist.* **51** 186–192.
[3] Newey, W. K. and Powell, J. L. (1990). Efficient estimation of linear and type I censored regression models under conditional quantile restrictions. *Econometric Theory* **6** 295–317. MR1085576
[4] Ruppert, D. (2002). Selecting the number of knots for penalized splines. *J. Comput. Graph. Statist.* **11** 735–757. MR1944261
[5] Sorva, R., Perheentupa, J. and Tolppanen, E.-M. (1984). A novel format for a growth chart. *Acta Paediatrica Scandinavica* **73** 527–529.
[6] Wei, Y., Pere, A., Koenker, R. and He, X. (2006). Quantile regression methods for reference growth charts. *Statistics in Medicine* **25** 1369–1382.

Department of Biostatistics
Columbia University
New York, New York 10032
USA
E-mail: yw2148@columbia.edu

Department of Statistics
University of Illinois
Champaign, Illinois 61820
USA
E-mail: x-he@uiuc.edu